\documentclass[12pt]{amsart}

\usepackage{environ, xcolor}
\NewEnviron{commentA}{}

\newcommand\Aoff{\RenewEnviron{commentA}{}}
\Aoff{} 

\usepackage{amsmath, amssymb}
\usepackage{array}
\usepackage[frame,cmtip,arrow,matrix,line,graph,curve]{xy}
\usepackage{graphpap, color, paralist, pstricks}
\usepackage[mathscr]{eucal}
\usepackage[pdftex]{graphicx}
\usepackage[pdftex,colorlinks,backref=page,citecolor=blue]{hyperref}
\usepackage{cleveref}
\usepackage{tikz, verbatim}

\newtheorem{theorem}{Theorem}[section]
\newtheorem{proposition}[theorem]{Proposition}
\newtheorem{corollary}[theorem]{Corollary}
\newtheorem{lemma}[theorem]{Lemma}

\theoremstyle{definition}
\newtheorem{definition}[theorem]{Definition}

\newtheorem{claim}[theorem]{Claim}

\newtheorem{question}[theorem]{Question}
\newtheorem{remark}[theorem]{Remark}

\usepackage[margin=1in]{geometry} 
\usepackage{amsmath,amssymb,mathtools, mathrsfs,esint,tikz-cd,verbatim}

\newcommand{\Z}{{\mathbb Z}}
\newcommand{\Q}{{\mathbb Q}}

\newcommand{\C}{{\mathbb C}}
\newcommand{\A}{{\mathbb A}}
\newcommand{\G}{{\mathbb G}}
\newcommand{\on}[1]{\operatorname{#1}}

\newcommand{\Spec}{{\on{Spec}}}

\newcommand{\set}[1]{\{#1\}}



\title{Quotient singularities in the Grothendieck ring of varieties}
\author{Louis Esser}
\author{Federico Scavia}

\address{Department of Mathematics, Princeton University, Fine Hall, Washington Road, Princeton, NJ 08544-1000, USA}
\email{esserl@math.princeton.edu}

\address{Institut Galil\'ee\\
	Universit\'e Sorbonne Paris Nord\\
	99 avenue Jean-Baptiste Cl\'ement, 93430\\ 
	Villetaneuse, France}
\email{scavia@math.univ-paris13.fr}

\begin{document}
	\maketitle
	
	\begin{abstract}
    Let $G$ be a finite group, $X$ be a smooth complex projective variety with a faithful $G$-action, and $Y$ be a resolution of singularities of $X/G$. Larsen and Lunts asked whether $[X/G]-[Y]$ is divisible by $[\mathbb{A}^1]$ in the Grothendieck ring of varieties. We show that the answer is negative if $BG$ is not stably rational and affirmative if $G$ is abelian. The case when $X=Z^n$ for some smooth projective variety $Z$ and $G=S_n$ acts by permutation of the factors is of particular interest. We make progress on it by showing that $[Z^n/S_n]-[Z\langle n\rangle / S_n]$ is divisible by $[\mathbb{A}^1]$, where $Z\langle n\rangle$ is Ulyanov's polydiagonal compactification of the $n$-th configuration space of $Z$.
	\end{abstract}
	
\section{Introduction}

For a field $k$, the Grothendieck ring of varieties $K_0(\on{Var}_k)$ is generated as an abelian group by isomorphism classes $[X]$ of separated $k$-schemes $X$ of finite type, subject to the relations $[X \setminus Y] = [X] - [Y]$, where $Y$ is a closed subscheme of $X$. The product of classes is induced by the fiber product of varieties over $k$.

The class $[X]$ of a variety $X$ in $K_0(\on{Var}_k)$ retains a great deal of geometric information about $X$.  This information is often expressed by {\it motivic measures}, i.e., ring homomorphisms $K_0(\on{Var}_k) \rightarrow A$.  When $k$ has characteristic zero, one important example is the homomorphism
\[\on{sb}\colon  K_0(\on{Var}_k) \rightarrow \Z[\on{SB}_k].\]

Here $\on{SB}_k$ is the monoid of stable birational equivalence classes of smooth projective varieties over $k$, and $\on{sb}$ is the homomorphism that sends a smooth projective connected variety to its stable birational equivalence class. By a result of Larsen and Lunts \cite{LLsr}, $\on{sb}$ induces an isomorphism $K_0(\on{Var}_k)/(\mathbb{L})\xrightarrow{\sim}\Z[\on{SB}_k]$, where  $\mathbb{L} = [\A_k^1]$. Therefore, it's natural to consider classes of varieties $X$ in the quotient $K_0(\on{Var}_k)/(\mathbb{L})$, but the interpretation of such classes is less clear when $X$ is not smooth.

Let $Z$ be a smooth projective variety. 
By definition, the $n$th {\it symmetric power} $\on{Sym}^n(Z)$ is the quotient of $Z^n$ by the action of the symmetric group $S_n$ which permutes the factors of $Z^n$.  When $\dim(Z) > 1$ and $n > 1$, $\on{Sym}^n(Z)$ is singular. The classes of symmetric powers in the Grothendieck ring are of particular interest because they appear in the {\it motivic zeta function} of $Z$:
$$\zeta_Z(t) \coloneqq \sum_{n = 0}^{\infty} [\on{Sym}^n(Z)] t^n.$$

When $Z$ is defined over a finite field, 
$\zeta_Z(t)$ specializes to the ordinary Weil zeta function under the point-counting measure.  However, unlike the Weil zeta function, the motivic zeta function $\zeta_Z(t)$ is not always rational \cite[Theorem 1.1]{LL}. Not much is known in general about the motivic zeta function or the classes of symmetric powers. In particular, the following basic question of Larsen and Lunts is open.

\begin{question}(Larsen and Lunts \cite[Question 6.7]{LL})
\label{mainq}
(i) Let $Z$ be a smooth projective complex variety, $n\geq 1$ be an integer, and $Y \rightarrow \on{Sym}^n(Z)$ be a resolution of singularities. Does the equality
\[[Y] = [\on{Sym}^n(Z)]\]
hold in the ring $K_0(\on{Var}_{\mathbb{C}})/(\mathbb{L})$?

(ii) Let $X$ be a smooth projective complex variety with the action of a finite group $G$. If $Y \rightarrow X/G$ is a resolution of singularities, do we have
\[[Y] = [X/G]\]
in $K_0(\on{Var}_{\mathbb{C}})/(\mathbb{L})$?
\end{question}

Question (i) is a special case of (ii). 

\Cref{mainq} arose as part of the study of the rationality of the motivic zeta function in \cite{LL}. Let $\Z[s]$ be a polynomial ring in one variable $s$ and $M\subset \Z[s]$ be the sub-monoid consisting of polynomials with constant term $1$. Larsen and Lunts considered a motivic measure \[\mu\colon K_0(\on{Var}_{\mathbb{C}})\to \Z[M]\] which sends the class of a $d$-dimensional smooth projective irreducible complex variety $X$ to $\sum_{i=0}^dh^{i,0}(X)s^i$. This homomorphism factors through $K_0(\on{Var}_k)/(\mathbb{L})$.  Both $K_0(\on{Var}_{\mathbb{C}})$ and $\Z[M]$ have natural structures as $\lambda$-rings, and it is natural to expect that $\mu$ is a $\lambda$-homomorphism. As observed by Larsen and Lunts, an affirmative answer to part (i) implies that $\mu$ is a $\lambda$-homomorphism. 

Our first main result shows that the answer to \Cref{mainq}(ii) is negative. Let $G$ be a finite group and let $k$ be a field. We say that the classifying space $BG$ is not stably rational over $k$ if the $k$-variety $V/G$ is not stably rational for a faithful $G$-representation $V$.  This definition is closely related to the famous Noether's problem on rationality of fields of invariants. It is independent of the choice of $V$ because the stable birational equivalence class of $V/G$ does not depend on the faithful representation by Bogomolov-Katsylo's double fibration method (see, e.g., the proof of \cite[Theorem 2.5]{totaro}). 

There are many known examples of groups $G$ with $BG$ not stably rational. Swan \cite{swan} and Voskresenski\u{\i} \cite{V} independently showed that $B(\Z/47\Z)$ is not stably rational over $\Q$. Lenstra \cite{lenstra} later classified the finite abelian groups $A$ such that $BA$ is stably rational over $k$; in particular, he showed that $B(\Z/8\Z)$ is not stably rational over $\Q$. When $k$ is algebraically closed, the first examples of groups with this property were constructed by Saltman \cite{saltman} and Bogomolov \cite{Bogomolov}. For example, if $n$ is any positive integer such that $2^6| n$ or $p^5|n$ for some odd prime $p$, then there exists a group $G$ of order $n$ such that $BG$ is not stably rational over $k$ \cite[Theorem 1.13]{HKK}. 

For any $G$ such that $BG$ is not stably rational over $k$, we produce an example answering \Cref{mainq}(ii) in the negative.

\begin{theorem}\label{counterexample}
     Let $k$ be a field of characteristic zero, $G$ be a finite $k$-group, $V$ be a faithful $k$-linear representation of $G$ containing at least one copy of the trivial representation, and $Y\to \mathbb{P}(V)/G$ be a resolution of singularities. Then $[Y] = [\mathbb{P}(V)/G]$ in $K_0(\on{Var}_k)/(\mathbb{L})$ if and only if $BG$ is stably rational over $k$. 
\end{theorem}

Conversely, we show that \Cref{mainq}(ii) has a positive answer when $G$ is abelian. This means that the motivic measure $\on{sb}$ can be used directly to detect the stable birational equivalence class of $X/G$ when $X$ is a smooth complex projective variety and $G$ is abelian. We refer the reader to \Cref{l-rational-def} for the definition of $\mathbb{L}$-rational singularities.

\begin{theorem}\label{thm-diag-group}
		Let  $G$ be a finite discrete abelian group of exponent $e$, $k$ be a field of characteristic zero containing a primitive $e$th root of unity, and $X$ be a smooth projective $k$-variety with a faithful $G$-action. Then $X/G$ has $\mathbb{L}$-rational singularities. In particular, if $Y$ is a resolution of singularities of $X/G$, then $[Y] = [X/G]$ in $K_0(\on{Var}_{k})/(\mathbb{L})$.
	\end{theorem}

Finally, we make progress toward answering \Cref{mainq}(i). For every $n\geq 1$, let $X\langle n \rangle$ be the polydiagonal compactification of the configuration space of $n$ distinct points on $X$, first constructed by Ulyanov \cite{Ulyanov}. This model, which is ``larger" than the Fulton-Macpherson compactification $X[n]$ of \cite{FM},  is obtained by successively blowing up certain diagonals in $X^n$ in an $S_n$-equivariant way. In particular, $X\langle n \rangle$ admits an $S_n$-action and the natural morphism $X\langle n \rangle\to X^n$ is $S_n$-equivariant. This map factors as $X \langle n \rangle \to X[n] \to X^n$. We refer the reader to \Cref{l-rational-def} for the definition of $\mathbb{L}$-rational fibers.

\begin{theorem}
\label{compactification-theorem}
Let $k$ be a field of characteristic zero. Then the morphism $X \langle n \rangle /S_n \rightarrow X^n/S_n$ has $\mathbb{L}$-rational fibers.  In particular, $[X\langle n\rangle/S_n]= [X^n/S_n]$ in $K_0(\on{Var}_k)/(\mathbb{L})$.
\end{theorem}

This does not answer \Cref{mainq}(i) because $X\langle n\rangle/S_n$ is typically singular. However, the stabilizers of the $S_n$-action on $X\langle n \rangle$ are simpler than on $X^n$, in fact even simpler than on $X[n]$: they are abelian \cite[Theorem 3.11]{Ulyanov}.  In particular, the singularities of $X \langle n \rangle/S_n$ may be amenable to study by toroidal methods.

\subsection*{Notation}
We work over a field $k$, with algebraic closure $\overline{k}$. A $k$-variety is a separated integral $k$-scheme of finite type. If $X$ is a $k$-scheme and $x$ is a point of $X$, we write $k(x)$ for the residue field of $x$. 

A $k$-group is a smooth affine group scheme of finite type over $k$. If $G$ is a finite $k$-group and $X$ is a quasi-projective $k$-scheme with a $G$-action, we denote by $X/G$ the quotient $k$-scheme.  We let $\G_{\on{m},X} \coloneqq \G_{\on{m},k}\times_k X$ denote the multiplicative group scheme over $X$.

\noindent{\it Acknowledgements. } The first author was partially supported by NSF grant DMS-2054553. Thank you to Daniel Litt and Burt Totaro for useful conversations.

\section{Preliminaries}

\subsection{Quotients in the Grothendieck Ring}

We first prove some properties of quotients that will be used throughout the paper. We also use these results to prove \Cref{counterexample}. We do not assume that $\on{char}(k)=0$ until it is necessary.
	
	\begin{lemma}\label{reduce-to-one-component}
		Let $k$ be a field, $G$ be a finite $k$-group, $X=\coprod_{i=1}^r X_i$ be a disjoint union of quasi-projective $k$-schemes such that $G$ acts on $X$ and the $G$-action transitively permutes the $X_i$. Let $G_1\subset G$ be the stabilizer of $X_1$. Then the open and closed embedding $X_1 \hookrightarrow X$ induces an isomorphism $X_1/G_1 \xrightarrow{\cong} X/G$.  Furthermore, if $f \colon Y\to X$ is a $G$-equivariant morphism, with $Y_i \coloneqq f^{-1}(X_i)$ and $f_i \coloneqq f|_{Y_i}$, then  $X_1\hookrightarrow X$ and $Y_1\hookrightarrow Y$ induce a commutative diagram
		\[
		\begin{tikzcd}
			Y_1/G_1 \arrow[r, "\cong"] \arrow[d]  & Y/G \arrow[d]\\
			X_1/G_1 \arrow[r, "\cong"] & X/G,
		\end{tikzcd}
		\]
		where the vertical maps are induced by $f_1$ and $f$, respectively.
	\end{lemma}
	
		\begin{proof}
		The commutativity of the diagram is clear, so it suffices to show that the horizontal maps are isomorphisms. For this, we are allowed to base change to $\bar{k}$, hence we may assume that $k$ is algebraically closed; in particular, $G$ is a finite discrete group. We may also suppose that $X$ and $Y$ are affine. The conclusion is implied by the following simple observation: if a finite group $G$ acts on a $k$-algebra $A$ via $k$-algebra automorphisms, $A=\prod_{i=1}^n A_i$,  the $A_i$ are transitively permuted by the $G$-action, and $G_1\subset G$ is the stabilizer of $A_1$, then $A^G = (A_1)^{G_1}$. 
	\end{proof}
	
	Let $G$ be a finite $k$-group, $X$ be a $k$-scheme with a $G$-action, and $\pi \colon P\to X$ be a $\G_{\on{m},X}$-torsor. We say that $\pi$ is a $G$-equivariant $\G_{\on{m},X}$-torsor if $G$ acts on $P$ and the $G$-action on $P$ commutes with the $\G_{\on{m},X}$-action, that is, the action map $\G_{\on{m},X}\times_X P\to P$ is $G$-equivariant, where $G$ acts on $\G_{\on{m},X}=\G_{\on{m},k}\times_k X$ via its action on $X$. By restricting the action map to $\set{1}\times P$, this implies that $\pi$ is  $G$-equivariant.
	
	\begin{proposition}\label{g-equiv-torsor}
		Let $k$ be a field, $G$ be a finite $k$-group, $X$ be a quasi-projective $k$-scheme of finite type, and $\pi \colon P\to X$ be a $G$-equivariant $\G_{\on{m},X}$-torsor. Then
		\[[P/G]=(\mathbb{L}-1)[X/G]\]
		in $K_0(\on{Var}_k)$.
	\end{proposition}
	
	\begin{proof}
	The conclusion of \Cref{g-equiv-torsor} is clear if $G$ is the trivial group because $\G_{\on{m},X}$-torsors are Zariski-locally trivial. When $G$ is arbitrary, by induction  we may suppose that the conclusion of \Cref{g-equiv-torsor} is true for all proper quotients of $G$. The proposition also holds for the group $G$ and torsors over quasi-projective schemes of dimension zero, since we can use \Cref{reduce-to-one-component} to reduce to a single point. Therefore, by induction, we can also assume that the conclusion of the proposition holds for the restriction of $\pi \colon P \rightarrow X$ to any $G$-invariant closed subscheme $Y \subset X$ of strictly smaller dimension. The scheme $X$ may consist of many irreducible components that intersect one another, but by the scissor relations in $K_0(\on{Var}_{k})$, we can remove the union of the intersections (which is $G$-invariant and of smaller dimension) and apply the inductive hypothesis. This reduces the problem to the case that $G$ transitively permutes the irreducible components of $X$, and then by \Cref{reduce-to-one-component} we may even assume that $X$ is irreducible.
		
		\begin{claim}\label{claim1}
			We may assume that $G$ acts faithfully on $P$.	
		\end{claim}
		
		Let $K\subset G$ be the kernel of the $G$-action on $P$. Then $K$ is contained in the kernel of the $G$-action on $X$, hence we may regard $\pi$ as a $G/K$-equivariant $\G_{\on{m},X}$-torsor. If $K$ is a non-trivial subgroup of $G$, then $G/K$ is a proper quotient of $G$ and the conclusion follows from the inductive assumption. This means that we may assume that $K$ is trivial, thus proving \Cref{claim1}.
		
		\begin{claim}\label{claim2}
			We may assume that all $\overline{k}$-points of $X$ have the same $G$-stabilizer $H$ and that $G$ acts freely on $P$.	
		\end{claim}
		
		Let $H\subset G$ be the kernel of the $G$-action on $X$, so that $H$ is a normal subgroup. Since $G$ is finite and $X$ is irreducible, there exists a dense open subscheme $U\subset X$ such that for all geometric points $x\in U(\overline{k})$, the $G$-stabilizer of $x$ is equal to $H$. The complement $X \setminus U$ is a union of subvarieties of lower dimension, hence we may replace $X$ by $U$ and $P$ by $P|_U$, so that the $G$-stabilizer of every $x\in X(\overline{k})$ is equal to $H$. In particular, every $G$-stabilizer in $P$ is a subgroup of $H$, and so coincides with the $H$-stabilizer. 
		
		By \Cref{claim1}, we may assume that $G$ acts faithfully on $P$. Since $G$ is finite, this implies the existence of a dense open subscheme $V$ of $P$ such that $G$ acts freely on $V$.
		
		Since $G$ acts on $P$ via $\G_{\on{m},X}$-torsor automorphisms and $H$ acts trivially on $X$, for all $x\in X(\overline{k})$, $p\in P_x(\overline{k})$, $\lambda\in \G_{\on{m},k}(\overline{k})$, and $h\in H$, we have $h\cdot p \in P_x(\overline{k})$ and
		\[h\cdot p=p\quad \Rightarrow\quad  h\cdot (\lambda p)=\lambda (h\cdot p)=\lambda p.\]
  
		 This shows that the $H$-stabilizer of $p$ is contained the $H$-stabilizer of $\lambda  p$. Replacing $p$ by $\lambda p$ and $\lambda$ by $\lambda^{-1}$ then shows that the $H$-stabilizer of $p$ is equal to the $H$-stabilizer of $\lambda p$. In other words, the $G$-stabilizers of geometric points of $P$ (which as we said in the previous paragraph coincide with their $H$-stabilizers) are constant along each fiber of $\pi$. In particular, $G$ acts freely on the open subscheme $\G_{\on{m},k}\cdot V=\pi^{-1}(\pi(V))\subset P$. By the inductive assumption, we are allowed to replace $X$ by the open subscheme $\pi(V)$ and $P$ by $P|_{\pi(V)}$. This proves \Cref{claim2}. 
		
		We are now ready to complete the proof of \Cref{g-equiv-torsor}. We may assume that \Cref{claim2} holds. 
		We may also replace $X$ by a dense open subscheme over which $\pi$ is trivial, that is, we may identify $\pi$ with the second projection $\G_{\on{m},X}\to X$.	 The $H$-action on the fibers of $\pi$ induces an $X$-group embedding $H_X\hookrightarrow \G_{\on{m},X}$. By the classification of isotrivial groups of multiplicative type \cite[Expos\'e X, Corollaire 1.2]{SGA3}, the proper isotrivial $X$-subgroups of $\G_{\on{m},X}$ are all of the form $\mu_{n,X}$ for some integer $n\geq 1$. Therefore, $H_X$ is $X$-isomorphic to $\mu_{n,X}$ for some integer $n\geq 1$, and this isomorphism identifies the action of $H_X$ on $P$ with the $\mu_{n,X}$-action obtained by restricting the $\G_{\on{m},X}$-torsor action. 
  (If $X(k)$ is empty, $H_X\cong \mu_{n,X}$ does not necessarily imply that $H\cong\mu_n$ over $k$.) 
		
		Therefore, the morphism $P/G\to X/G$ can be obtained by first taking the quotient of $P$ by the action of  $H_X \cong \mu_{n,X}$, and then taking the quotient by the action of  $(G/H)_X$. Said otherwise, $P/G\to X/G$ is the same as $(P/\mu_{n,X})/(G/H)\to X/(G/H)$, that is, the quotient of $P/\mu_{n,X}\to X$ by the induced $G/H$-action. Note that the $\G_{\on{m},X}$-action of scalar multiplication on $P$ descends to a $\G_{\on{m},X}$-action on $P/\mu_{n,X}\to X$, with action map fitting into the commutative square
		\[
		\begin{tikzcd}
			\G_{\on{m},X}\times_XP\arrow[r] \arrow[d]  & P \arrow[d]  \\
			\G_{\on{m},X}\times_X(P/\mu_{n,X})\arrow[r] & P/\mu_{n,X}, 	
		\end{tikzcd}
		\]
		where the vertical map on the left sends $(t,p)\mapsto (t^n, [p])$, and the map on the right sends $p\mapsto [p]$. (Equivalently, we quotient both $\G_{\on{m},X}$ and $P$ by $\mu_{n,X}$ and then identify $\G_{\on{m},X}/\mu_{n,X}$ with $\G_{\on{m},X}$.) 
		
		If $H$ is a non-trivial subgroup of $G$, then the conclusion follows from the inductive assumption. If $H$ is the trivial subgroup of $G$, then by \Cref{claim2} the group $G$ acts freely on $X$ and $P$. In this case, by descent along \'etale-torsors \cite[Theorem 4.46]{vistoli-descent} the quotient $P/G\to X/G$ is a $\G_{\on{m},X/G}$-torsor, hence a Zariski-locally trivial fibration with fibers isomorphic to $\A^1 \setminus \set{0}$, and we conclude that $[P/G]=(\mathbb{L}-1)[X/G]$.
	\end{proof}

	\begin{corollary}\label{cor-tautological-bundle}
		Let $G$ be a finite $k$-group, $W$ be a finite $k$-scheme with a $G$-action such that $W/G=\on{Spec}(k)$, and $V$ be a $G$-equivariant vector bundle on $W$. Then
		\[[V/G]=(\mathbb{L}-1)[\mathbb{P}(V)/G]+1\]
		in $K_0(\on{Var}_k)$. In particular, the following holds in $K_0(\on{Var}_k)/(\mathbb{L})$: $[\mathbb{P}(V)/G] = 1$ if and only if $[V/G] = 0$.
	\end{corollary}
	
	\begin{proof}
	Consider the $\G_{\on{m}}$-torsor $\pi \colon V \setminus s(W) \to \mathbb{P}(V)$, where $s \colon  W \rightarrow V$ is the zero section of the vector bundle. Since $G$ acts linearly on $V$ over $W$, the $G$-action commutes with the $\G_{\on{m}}$-action, that is, $\pi$ is a $G$-equivariant $\G_{\on{m}}$-torsor. Therefore, the assumptions of \Cref{g-equiv-torsor} are satisfied. By \Cref{g-equiv-torsor}, we have
		\[[V/G]=[(V \setminus s(W))/G]+[s(W)/G]=(\mathbb{L}-1)[\mathbb{P}(V)/G]+1.\qedhere\]
	\end{proof}
	
	\Cref{cor-tautological-bundle} is an important step towards our proof of \Cref{compactification-theorem}. It also allows us to prove \Cref{counterexample}.

	\begin{proof}[Proof of \Cref{counterexample}]
    By assumption, $V=V'\oplus k$ for some $k$-linear representation $V'$ of $G$. Note that $V'$ is faithful because $V$ is faithful and $G$ acts trivially on $k$. The variety $\mathbb{P}(V)$ is a union of the $G$-invariant affine chart isomorphic to $V'$ and its complement $\mathbb{P}(V')$, both with the natural $G$-action.  Therefore, in $K_0(\on{Var}_k)$ we have
    \[[\mathbb{P}(V)/G] = [V'/G] + [\mathbb{P}(V')/G] = (\mathbb{L}-1)[\mathbb{P}(V')/G] + 1 + [\mathbb{P}(V')/G] = \mathbb{L}[\mathbb{P}(V')/G] + 1,\]
    where the second equality follows from \Cref{cor-tautological-bundle}. This implies that $[\mathbb{P}(V)/G] = 1$ in $K_0(\on{Var}_k)/(\mathbb{L})$. However, $\mathbb{P}(V)/G$ contains the open subvariety $V'/G$ which is stably rational if and only if $BG$ is. Therefore, if $Y \rightarrow \mathbb{P}(V)/G$ is any resolution, $[Y] \equiv 1$ in $K_0(\on{Var}_k)/(\mathbb{L})$ if and only if $Y$ is stably rational, because $Y$ is a smooth projective variety \cite{LLsr}.
	\end{proof}
	
	We conclude this section with the observation that \Cref{g-equiv-torsor} allows one to compute the class of the quotient of a $G$-equivariant line bundle.
	
	\begin{remark}\label{equiv-line-bundle}
    Let $k$ be a field, $G$ be a finite group, $X$ be a quasi-projective $k$-scheme of finite type, and $L\to X$ be a $G$-equivariant line bundle. Let $s \colon X\to L$ be the zero section of $L$, and let $P \coloneqq L \setminus s(X)$. Then $L/G$ is the union of a closed subscheme isomorphic to $X/G$ and with complement $P/G$. (This is automatic when $\on{char}(k)=0$, but is true in arbitrary characteristic: the $G$-equivariant surjection $\mathcal{O}_P\to \mathcal{O}_X$ induced by $s$ is $G$-equivariantly split, and so the induced map $\mathcal{O}_P^G\to \mathcal{O}_X^G$ is also surjective.) Since $P$ is a $G$-equivariant $\G_{\on{m},X}$-torsor,  \Cref{g-equiv-torsor} implies
    \[[L/G]=[P/G]+[X/G]=\mathbb{L}[X/G]\]
    in $K_0(\on{Var}_k)$.
	\end{remark}

\subsection{Toroidal embeddings and \texorpdfstring{$\mathbb{L}$}{L}-rational singularities}

In this section, we'll introduce $\mathbb{L}$-rational singularities and some of the language of toroidal embeddings. These notions will be used to prove the equalities of classes in $K_0(\on{Var}_k)/(\mathbb{L})$ given in \Cref{thm-diag-group} and \Cref{compactification-theorem}.
	
\begin{definition}\label{l-rational-def}{(cf. \cite[Definition 4.2.4]{NS})}
		Let $k$ be a field and $X$ and $Y$ be $k$-varieties.  We say that a proper morphism of schemes $h \colon  Y \rightarrow X$ has {\it $\mathbb{L}$-rational fibers} if for each point $x$ of the scheme $X$, $[h^{-1}(x)] = 1$ in $K_0(\on{Var}_{k(x)})/(\mathbb{L})$.  We say that $X$ has {\it $\mathbb{L}$-rational singularities} if there exists a resolution of singularities $h \colon  Y \rightarrow X$ with $\mathbb{L}$-rational fibers.
	\end{definition}
	
	These definitions are local on $X$ in the Zariski topology. As observed after \cite[Definition 4.2.4]{NS}, the weak factorization theorem in characteristic zero \cite{AKMW} implies that the definition of $\mathbb{L}$-rational singularities is independent of the choice of resolution of $X$. The following lemma relates this local definition to the global equality of classes:
	
\begin{lemma}{\cite[Lemma 4.2.5]{NS}}
\label{spreadingout}
Let $k$ be a field and $h \colon  Y \rightarrow X$ be a proper morphism of $k$-varieties.  If $h$ has $\mathbb{L}$-rational fibers, then $[Y] = [X]$ in $K_0(\on{Var}_k)/(\mathbb{L})$.
\end{lemma}
	
In particular, this implies that if $X$ has $\mathbb{L}$-rational singularities, the class of $X$ equals the class of its resolution modulo $\mathbb{L}$.

\begin{lemma}
\label{closure}
The class of morphisms with $\mathbb{L}$-rational fibers is closed under base change and composition.
\end{lemma}

\begin{proof}
First, we prove stability under base change.  Suppose that $h \colon  X \rightarrow Y$ is a proper morphism with $\mathbb{L}$-rational fibers and that $g \colon  S \rightarrow Y$ is an arbitrary morphism of $k$-varieties.  Let $f \colon  X \times_Y S \rightarrow S$ be the base change, $s \in S$ be a point of the scheme $S$, and $y = g(s)$ the image in $Y$.  We have a commutative diagram

\begin{center}
\begin{tikzcd}
    f^{-1}(s) \arrow[r] \arrow[d] & X \times_Y S \arrow[r] \arrow[d,"f"] & X \arrow[d, "h"] \\
    \Spec(k(s)) \arrow[r] & S \arrow[r, "g"] & Y.
\end{tikzcd}
\end{center}

Since the two inner squares are pullbacks, the outer rectangle is as well.  The morphism $\Spec(k(s)) \rightarrow Y$ factors as $\Spec(k(s)) \rightarrow \Spec(k(y)) \rightarrow Y$, where $k(s)/k(y)$ is an extension of fields.  Thus, in the new diagram
\begin{center}
\begin{tikzcd}
     f^{-1}(s) \arrow[r] \arrow[d] & h^{-1}(y) \arrow[r] \arrow[d,"f"] & X \arrow[d] \\
    \Spec(k(s)) \arrow[r] & \Spec(k(y)) \arrow[r, "g"] & Y,
\end{tikzcd}
\end{center}
the outer and right rectangles are pullbacks, so the left is too.  Since $[h^{-1}(y)] = 1$ in $K_0(\on{Var}_{k(y)})/(\mathbb{L})$ by assumption, we have $[f^{-1}(s)] = 1$ in $K_0(\on{Var}_{k(s)})/(\mathbb{L})$ by functoriality of the Grothendieck ring.

Now we'll prove stability under composition. Suppose that $g \colon  X \rightarrow Y$ and $h \colon  Y \rightarrow Z$ are proper with $\mathbb{L}$-rational fibers and let $z \in Z$ be a point.  We have $[h^{-1}(z)] = 1$ in $K_0(\on{Var}_{k(z)})/(\mathbb{L})$ by assumption.  Further, $g^{-1}(h^{-1}(z)) \rightarrow h^{-1}(z)$ has $\mathbb{L}$-rational fibers by stability under base change.  Then by \Cref{spreadingout}, $[g^{-1}(h^{-1}(z))] = [h^{-1}(z)] = 1$ in $K_0(\on{Var}_{k(z)})/(\mathbb{L})$, completing the proof.
\end{proof}
	
As a consequence, we also obtain some useful properties of $\mathbb{L}$-rational singularities.

\begin{lemma}
\label{singcompclosure}
Let $g \colon  X' \rightarrow X$ be a proper birational morphism of $k$-varieties with $\mathbb{L}$-rational fibers.  Then, if $X'$ has $\mathbb{L}$-rational singularities, so does $X$.	
\end{lemma}

\begin{proof}
Let $h \colon  Y \rightarrow X'$ be a resolution of singularities, which has $\mathbb{L}$-rational fibers by assumption.  Since $g$ is birational, $g \circ h \colon  Y \rightarrow X$ is also a resolution of singularities.  The morphism $g \circ h$ is also proper with $\mathbb{L}$-rational fibers by \Cref{closure}, completing the proof.
\end{proof}

\begin{lemma}
\label{etalebasechange}
If a $k$-variety $X$ has $\mathbb{L}$-rational singularities and $f \colon  X' \rightarrow X$ is an \'{e}tale morphism, then $X'$ has $\mathbb{L}$-rational singularities as well.
\end{lemma}

\begin{proof}
    If $Y \rightarrow X$ is a resolution of singularities, then the pullback $Y' = Y \times_X X' \rightarrow X'$ is also a resolution of singularities. This is because the properties of being proper and birational are preserved under pullback by an \'{e}tale morphism, while $Y'$ is an \'{e}tale cover of a smooth variety $Y$ and hence is smooth.  By \Cref{closure}, $Y' \rightarrow X'$ has $\mathbb{L}$-rational fibers, completing the proof.
\end{proof}

Next, we'll introduce some of the language of toroidal embeddings.  This machinery is used for the proof of \Cref{thm-diag-group} in \Cref{section-diag-groups}, but is not required for \Cref{compactification-theorem}. The datum of a {\it toroidal embedding} consists of a normal variety $X$ over a field $k$ and a non-empty open subscheme $U \subset X$ with the property that $(X,U)$ looks \'{e}tale-locally like the embedding of the open torus orbit in a toric variety.  That is, for every point $x$ of the scheme $X$, there is an \'{e}tale neighborhood $f \colon  V \rightarrow X$ of $x$ and a toric variety $(X_{\sigma},T_{\sigma})$ such that $V$ admits an \'{e}tale morphism $g \colon  V \rightarrow X_{\sigma}$ and $g^{-1}(T_{\sigma}) = f^{-1}(U)$.  A toroidal embedding is {\it strict} if we instead require that for each $x \in X$ there is a Zariski open neighborhood $V$ of $x$ satisfying the same property.
	
There are equivalent definitions of toroidal and strictly toroidal originally due to Mumford (\cite[p. 195, footnote]{KKMSD}; see also \cite[Section 2]{ADK} or \cite[Section 4]{W} for more on these definitions).
	
	\begin{lemma}{\cite[p. 195, footnote]{KKMSD}, \cite[Lemma 4.8.3]{W}}
	\label{toroidal-def}
	    Let $k$ be a field of characteristic zero, $\bar{k}$ be an algebraic closure of $k$, $X$ be a normal $k$-variety, and $U$ be a non-empty open subscheme of $X$.
	    \begin{enumerate}
	        \item   The pair $(X,U)$ is a toroidal embedding if and only if for every closed point $x$ in $X_{\bar{k}}$, there exists a toric variety $(X_{\sigma},T_{\sigma})$ over $\bar{k}$, a point $t \in X_{\sigma}$ and an isomorphism 
	        $$\varphi \colon  \widehat{\mathcal{O}}_{X_{\bar{k}},x} \rightarrow \widehat{\mathcal{O}}_{X_{\sigma},t}$$
	        of complete local algebras over $\bar{k}$ with the property that $\varphi$ takes the ideal of $X_{\bar{k}} \setminus U_{\bar{k}}$ to the ideal of $X_{\sigma} \setminus T_{\sigma}$.
	        \item A toroidal embedding $(X,U)$ is strict if and only if each irreducible component of $X \setminus U$ is normal.
	    \end{enumerate}
	\end{lemma}
 
	Strict toroidal embeddings enjoy many of the properties of toric varieties. In particular, the singularities of strictly toroidal embeddings are $\mathbb{L}$-rational \cite[Example 4.2.6(3)]{NS}.  This is because toric singularities are $\mathbb{L}$-rational and strictly toroidal embeddings Zariski-locally admit \'{e}tale morphisms to toric varieties, so we may apply \Cref{etalebasechange}.
	
	Let $(X,U)$ be a toroidal embedding, and let $G$ be a finite discrete group acting on $X$ and leaving $U$ invariant. Following \cite[Section 2.3]{ADK}, we say that the action of $G$ on $(X,U)$ is {\it toroidal} if for each closed point $x$ in $X_{\bar{k}}$ the following holds: there exists a local model $(X_{\sigma},T_{\sigma})$ as in \Cref{toroidal-def}(1) and a homomorphism $\alpha \colon  G_x \rightarrow T_{\sigma}$ from the stabilizer $G_x$ of $x$ under the induced action to the torus such that the action of $G_x$ on $\widehat{\mathcal{O}}_{X_{\sigma},t}$ via the isomorphism $\varphi$ agrees with the action on the toric variety $X_{\sigma}$ via $\alpha$. We say the $G$-action on $(X,U)$ is {\it strict} if $X \setminus U = D$ has the property that $\bigcup_{g \in G} g \cdot Z$ is normal for every irreducible component $Z$ of $D$.  Equivalently, $\bigcup_{g \in G} g \cdot Z$ is a disjoint union of normal components, so that either $g \cdot Z = Z$ or $g \cdot Z \cap Z = \emptyset$. 
	
	In the event that we have a strict toroidal action on a strict toroidal embedding, the quotient also inherits a strict toroidal structure (see \cite[Section 2.3]{ADK}, \cite[Proposition 2.5]{AW}):
	
	\begin{proposition}
		\label{quotient}
		If $(X,U)$ is a strict toroidal embedding on which $G$ acts strictly and toroidally, then $(X/G,U/G)$ is also a strict toroidal embedding.
	\end{proposition}

 \begin{commentA}

 Here is the proof of \Cref{quotient}.
	\begin{proof}
	  To see that $(X/G,U/G)$ is toroidal, it's enough to consider the situation on an \'{e}tale-local toric model over the base change to $\bar{k}$, where the result follows from the fact that the quotient of a toric variety by a finite subgroup of the torus is still toric.  
	  
	  To prove that the toroidal embedding $(X/G,U/G)$ is strict, it suffices by \Cref{toroidal-def}(2) to show that any irreducible component of the complement $(X/G) \setminus (U/G)$ is normal.  Since the $G$-action is strict, any component of the complement  is of the form $(\bigcup_{g \in G} g\cdot Z)/G$ for some irreducible component $Z \subset X \setminus U$ and this is a union of disjoint normal varieties (possibly with repeats).  By \Cref{reduce-to-one-component}, the quotient is $Z/H$ where $H$ is the stabilizer of $Z$, hence it is normal.
	\end{proof}
\end{commentA}
	
	\section{Diagonalizable Groups}
    \label{section-diag-groups}

    Recall that a $k$-group $G$ is {\it diagonalizable} if it embeds in $\G_{\on{m},k}^n$ for some $n\geq 1$. If $G$ is an abstract finite group, then $G$ is diagonalizable if and only if $G$ is abelian and $k$ contains a primitive root of unity of order equal to the exponent of $G$.
    
    Let $X$ be a smooth projective variety.  In this section, we prove \Cref{thm-diag-group}, namely that the singularities of a quotient $X/G$ are $\mathbb{L}$-rational in the special case that the abstract finite group $G$ is diagonalizable. As a first step, we see that an equivariant blowup of $X$ in a smooth subvariety leaves the class of the quotient unchanged.  This result holds for any (not necessarily discrete) diagonalizable $k$-group $G$.
	
	\begin{proposition}
		\label{blowup}
		Suppose that $G$ is a diagonalizable $k$-group acting on a smooth quasi-projective variety $X$ and $f \colon  Y \rightarrow X$ is an equivariant blowup of $X$ in a smooth proper subvariety $Z$.  Then the induced morphism $h \colon  Y/G \rightarrow X/G$ has $\mathbb{L}$-rational fibers.  In particular,
		$$[Y/G] = [X/G]$$
		holds in $K_0(\on{Var}_{k})/(\mathbb{L})$.
	\end{proposition}
	
	\begin{proof}
	Since $f$ is equivariant, we must have that the action of $G$ on $X$ restricts to an action of $G$ on $Z$. Over $Z$, the morphism $f$ is a projective bundle $\mathbb{P} (\mathcal{N}_{Z/X})$, where $\mathcal{N}_{Z/X}$ is the normal bundle to $Z$ in $X$. Denote by $p \colon  X \rightarrow X/G$ the quotient morphism.  We'll show that $[h^{-1}(p(x))] = 1$ in $K_0(\on{Var}_{k(p(x))})/(\mathbb{L})$ for any point $x \in X$. The quotient map $h \colon  Y/G \rightarrow X/G$ is an isomorphism away from $Z/G$, so we need only consider the case that $x \in Z$.  Suppose that $W$ is the orbit of the point $x$ under the $G$ action.  It is a finite disjoint union of points.
	
	The quotient $W/G$ is the point $\Spec(k(p(x)))$ and the fiber of $h$ over $p(x)$ is the quotient of the projectivization $\mathbb{P}N_W$ by $G$, where $N_W$ is the pullback of the normal bundle to the scheme $W$ (here we use that projectivization commutes with base change).  Therefore, we may apply \Cref{cor-tautological-bundle} to conclude that $[h^{-1}(p(x))] = [\mathbb{P}(N_W)/G] = 1$ if and only if $[N_W/G] = 0$ in the ring $K_0(\on{Var}_{k(p(x))})/(\mathbb{L})$. The total space $N_W$ is a disjoint union of vector spaces over each point in $W$, which are transitively permuted by the $G$-action.  Therefore, by \Cref{reduce-to-one-component}, $N_W/G \cong N_x/G_1$, where $G_1 \subset G$ is the subgroup that stabilizes the point $x \in W$.  It will suffice to prove that the class $[N_x/G_1]$ is trivial in $K_0(\on{Var}_{k(p(x))})/(\mathbb{L})$. 
	
	The fiber $N_x$ is a finite-dimensional vector space over $k(x)$, but $G_1$ may not act by $k(x)$-automorphisms because it also acts on the field $k(x)$ non-trivially.  However, $k(x)$ is a finite extension of the subfield $k(p(x)) = k(x)^{G_1}$ and $G_1$ does act by $k(p(x))$-automorphisms.  Since $G_1 \subset G$ is diagonalizable over $k$, it is certainly diagonalizable over $k(p(x))$, so $N_x/G_1$ is the quotient of a finite-dimensional vector space by a diagonalizable group action.  This has trivial class in $K_0(\on{Var}_{k(p(x))})/(\mathbb{L})$ by the following lemma.
	
	\begin{lemma}
	    Let $V$ be a vector space over a field $k$ with a diagonalizable action by a $k$-group $H$.  Then $[V/H] = 0$ in $K_0(\on{Var}_k)/(\mathbb{L})$.
	\end{lemma}
	
	\begin{proof}
		Using the scissor relations, we can decompose $V$ into coordinate strata, each of which is a split torus $(\G_{\on{m},k})^d$ with a diagonal action of $H$.  The quotient of this stratum by this action is again isomorphic to $(\G_{\on{m},k})^d$, so the class in the Grothendieck ring is unchanged.  Reassembling the pieces, we have $[V/H] = [V] = 0$ in $K_0(\on{Var}_{k})/(\mathbb{L})$.
	\end{proof}
	
	This shows that the fiber of $h \colon  Y/G \rightarrow X/G$ over the point $p(x)$ has class $1$ in $K_0(\on{Var}_{k(p(x))})/(\mathbb{L})$ for any point $x$ of the scheme $X$.  Since the quotient $X \rightarrow X/G$ is surjective, the image $p(x)$ ranges over every point of the scheme $X/G$ and the proof is complete.
	\end{proof}

	The proof of \Cref{thm-diag-group} relies on the existence of a $G$-equivariant blowup of $X$ on which the action of $G$ is well-behaved. We'll apply the main theorem of \cite{AW}:
	
	\begin{theorem}{\textup{(cf. \cite[Theorem 0.1]{AW})}}
		\label{res}
		Let $X$ be a projective variety over a field $k$ of characteristic zero with an action by a finite (discrete) group $G$. Then there is a $G$-equivariant proper birational morphism $r \colon  X_1 \rightarrow X$ and an open set $U \subset X_1$ such that $X_1$ is a nonsingular projective variety, $(X_1,U)$ has the structure of a strictly toroidal embedding, and $G$ acts strictly and toroidally on $(X_1,U)$.
	\end{theorem}
	
	The statement of the theorem differs slightly here from \cite{AW}, but this is what their proof shows.  That paper also assumes that $k$ is algebraically closed, but this assumption can be removed. 

\begin{commentA}
    For clarity, we reproduce some proof details below, following the same steps as in \cite{AW}.
	
	\begin{proof}
	Given $X$, $G$ as in the theorem, let $Y = X/G$ and $W$ be the branch locus of the quotient morphism.  Find a resolution of singularities $(Y',W') \rightarrow (Y,W)$ with the property that $W'$ is a strict normal crossings divisor.  Now let $X'$ be the normalization of $Y'$ in the function field $k(X)$ and $Z'$ the inverse image of $W'$.  Defining $U = X' \setminus Z'$, we have that $(X',U)$ is a strictly toroidal embedding on which the action of $G$ is strict toroidal (this follows from Abhyankar's lemma as in \cite[Lemma 3.3]{ADK}, with strictness coming from the fact that $W'$ has strict normal crossings).
	
	Now that we've produced a strictly toroidal model $(X',U)$, the last step is an equivariant toroidal resolution of singularities $X_1 \rightarrow X$, the existence of which is proven in \cite[Theorem 0.2]{AW}.  The proof works by finding a $G$-equivariant subdivision of the polyhedral complex $\Sigma$ associated to $(X',U)$ with nice properties.  This uses the correspondence between toroidal modifications of $(X',U)$ and subdivisions of the complex $\Sigma$, which was originally only proven for $k$ algebraically closed in \cite{KKMSD}.  However, the correspondence in fact holds over any field of characteristic zero \cite[Theorems 4.8.11, 4.11.2]{W}.
	\end{proof}
 \end{commentA}
	
	\begin{proof}[Proof of \Cref{thm-diag-group}] Let $G$ be any discrete, diagonalizable $k$-group.  Given a smooth $G$-variety $X$, find a birational morphism $X_1 \rightarrow X$ satisfying the conditions of \Cref{res}.  Then $X_1/G$ is strictly toroidal by \Cref{quotient}.  If $Y \rightarrow X_1/G$ is a resolution of singularities, we therefore have $[Y] = [X_1/G]$ in $K_0(\on{Var}_{\C})/(\mathbb{L})$.  But $Y$ is also a resolution of $X/G$.  By the $G$-equivariant weak factorization theorem, $X_1$ can be constructed from $X$ via some series of equivariant blowups and blowdowns in smooth centers.  Using \Cref{blowup} at each step, we have $[X_1/G] = [X/G]$ in $K_0(\on{Var}_{\C})/(\mathbb{L})$.  This completes the proof.
	\end{proof}

	\section{Symmetric Powers}
	\label{section-sym-powers}

    In this section, we'll prove \Cref{compactification-theorem}, which shows that the class of a symmetric power $\on{Sym}^n(X)$ in $K_0(\on{Var}_k)/(\mathbb{L})$ equals the class of a certain birational modification with simpler singularities.  This modification is the quotient of the polydiagonal compactification $X\langle n \rangle$ of the configuration space of $n$ distinct points on $X$ by $S_n$.  In \Cref{subsection-sym-reps}, we'll prove some general facts about quotients of symmetric group representations in the Grothendieck ring.  \Cref{subsection-polydiag} will introduce the polydiagonal compactification and complete the proof of \Cref{compactification-theorem}.

	\subsection{Stratification of Symmetric Group Representations}
    \label{subsection-sym-reps}
	
	We'll introduce some notation for partitions that will be used throughout the remainder of this section. For a positive integer $n$, let $P(n)$ denote the set of partitions of the set $[n] = \{1,\ldots,n\}$.  Suppose $\pi$ is a partition in $P(n)$. We say that $\pi$ is of {\it type} $\mathbf{a} = (1^{m_1},\dots,n^{m_n})$ if it contains exactly $m_h$ blocks of size $h$, $h = 1,\ldots,n$.  We'll write $\pi_1,\dots,\pi_i$ for the blocks of $\pi$, where $i = m_1 + \cdots + m_n$ is the total number of blocks.  For each block $\pi_j$, $|\pi_j|$ denotes the number of elements in the block.
	
	The group $S_n$ naturally acts on $[n] = \set{1,\dots,n}$, hence on the set of partitions $P(n)$. The orbits of this action are precisely the sets of partitions of a fixed type $\mathbf{a}$. Let $S_{\pi}\subset S_n$ be the stabilizer of $\pi$. We can naturally write $S_{\pi}$ as a semi-direct product
	\[S_{\pi}\cong S'_{\pi} \rtimes \overline{S}_{\pi}.\]
	
    Here $S'_{\pi} \cong S_{\pi_1} \times \cdots \times S_{\pi_{i}}$, where $S_{\pi_{l}}$ is the symmetric group on the elements of $\pi_{l}$.  The elements of the subgroup $S'_{\pi} \subset S_n$ therefore preserve each block individually.  The group $\overline{S}_{\pi}$ is defined as $S_{m_1} \times \cdots \times S_{m_n}$, where $S_{m_h}$ permutes the $m_h$ factors $S_{\pi_l}$ such that $|\pi_l|=m_h$.

    \begin{definition}
        We say that a morphism of $k$-varieties $Y\to Z$ is {\it stratified by vector bundles} (with height $N$) if it can be written as a composition \[Y=Y_N\to Y_{N-1}\to\dots \to Y_1\to Y_0=Z\] such that for all $i=1,\dots, N$ there exists a locally closed stratification of $Y_{i-1}$ where the restriction of $Y_i\to Y_{i-1}$ to each stratum is a vector bundle.
    \end{definition}
    
    It follows from the definition that the composition of two morphisms stratified by vector bundles is also stratified by vector bundles.  If a morphism $f \colon  Y \rightarrow Z$ is stratified by vector bundles, then $[f^{-1}(z)] = 0$ in $K_0(\mathrm{Var}_{k(z)})/(\mathbb{L})$ for every point $z$ in $Z$.
    
    Now let $X$ be a $k$-variety, $E\to X$ be a vector bundle of rank $r$, and $d$ be a positive integer. We denote by $(E^{\oplus d})_0\to X$ the $S_d$-equivariant vector subbundle of $E^{\oplus d}$ given by the kernel of the addition map $E^{\oplus d}\to E$. The symmetric group $S_d$ acts on $E^{\oplus d}$ over $X$ by permutation of the $d$ summands $E$, and this action leaves $(E^{\oplus d})_0$ invariant.
	
	\begin{lemma}\label{zero-sum-over-k}
		Let $r,d$ be non-negative integers.  Suppose that $S_d$ acts on $(\A^r)^d$ by permutation of the $d$ factors $\A^r$, and consider the induced $S_d$-action on $(\A^r)_0^d$.  The projection $p \colon  ((\A^r)^d_0)/S_d\to ((\A^{r-1})^d_0)/S_d$ is stratified by vector bundles with height $1$.
	\end{lemma}
	
	\begin{proof}
		We denote by $\iota_r \colon (\A^r)^d_0\hookrightarrow (\A^r)^d$ the canonical linear inclusion as the kernel of the addition map. We view the trivial vector bundle $(\A^r)^d$ as the space of $r$ by $d$ matrices
		\[
		A = \begin{bmatrix} 
			a_{1,1} & \dots  & a_{1,d}\\
			\vdots & \ddots & \vdots\\
			a_{r,1} & \dots  & a_{r,d} 
		\end{bmatrix}.
		\]
		In this notation, the $S_d$-action on $(\A^r)^d$ permutes the columns of $A$, the map $p$ is induced by forgetting the last row of $A$, and $(\A^r)^d_0$ parametrizes the matrices $A$ such that $a_{i,1}+\dots +a_{i,d}=0$ for all $1\leq i\leq r$.
		We obtain a commutative diagram of $S_d$-equivariant linear maps
		\begin{equation}
		\label{reduction}
		\begin{tikzcd}
			(\A^r)_0^d\arrow[dr, "p_0"] \arrow[r, hook]  \arrow[rr, bend left=30, hook, "\iota_r"] & \iota_{r-1}^*((\A^r)^d) \arrow[r, hook] \arrow[d,"\iota_{r-1}^*p"]   & (\A^r)^d \arrow[d, "p"] 	\\
			& (\A^{r-1})_0^d \arrow[r, hook, "\iota_{r-1}"] & (\A^{r-1})^d,
		\end{tikzcd}
		\end{equation}
        where the square on the right is cartesian. Note that $\iota_{r-1}^*((\A^r)^d)$ parametrizes those matrices $A$ for which $a_{i,1}+\dots +a_{i,d}=0$ for all $1\leq i\leq r-1$. The inclusion of $(\A^r)_0^d$ in $\iota_{r-1}^*((\A^r)^d)$ is defined by the extra condition $a_{r,1}+\dots +a_{r,d}=0$.
		
		Consider the $S_d$-equivariant locally-closed stratification $\set{Y_{\pi}}_{\pi\in P(d)}$ of $(\A^{r-1})^d$, where $Y_{\pi}$ parametrizes $d$-tuples $(v_1,\dots,v_d)\in (\A^{r-1})^d$ such that $v_s=v_t$ if and only if $s$ and $t$ are in the same block of $\pi$. The action of $S_d$ on the set of partitions $P(d)$ breaks into orbits by partition type. By \Cref{reduce-to-one-component}, $p$ is stratified by the union of maps of the form $p^{-1}(Y_{\pi})/S_{\pi}\to Y_{\pi}/S_{\pi}$, that is,
		\[(Y_{\pi}\times (\A^1)^d)/S_{\pi}\to Y_{\pi}/S_{\pi}.\]
		
		Here there is one stratum for each partition type $\mathbf{a} = (1^{m_1},\ldots,n^{m_n})$ and $\pi$ can be taken to be any representative of this partition type.
		
		We rearrange the factors of $(\A^1)^d$ by putting together the factors $\A^1$ corresponding to elements of $\set{1,\dots,d}$ in the same block of $\pi$. The previous map becomes
		\[(Y_{\pi}\times \prod_{j=1}^i \A^{|\pi_j|})/S_{\pi}\to Y_{\pi}/S_{\pi}.\]
  
		The subgroup $S'_{\pi}$ acts trivially on $Y_{\pi}$, and the direct factor $S_{\pi_j}\subset S'_{\pi}$ acts on $(\A^1)^{|\pi_j|}$ by permuting the coordinates and trivially on the other factors. By the fundamental theorem on symmetric polynomials, the quotient $(\A^1)^{|\pi_j|}/S_{\pi_j}$ is isomorphic to $\A^{|\pi_j|}$, and the quotient map $(\A^1)^{|\pi_j|}\to\A^{|\pi_j|}$ is given by the elementary symmetric functions. We have obtained a commutative diagram
		\begin{equation}\label{descent}
			\begin{tikzcd}
				Y_{\pi}\times (\A^1)^d\arrow[r, "\sigma"]\arrow[dr, "p"]  &  Y_{\pi}\times \prod_{j=1}^i \A^{|\pi_j|} \arrow[r] \arrow[d] & p^{-1}(Y_{\pi})/S_{\pi}=\overline{p}^{-1}(Y_{\pi}/S_{\pi}) \arrow[d, "\overline{p}"]  \\
				& Y_{\pi} \arrow[r] & Y_{\pi}/S_{\pi},	
			\end{tikzcd}
		\end{equation}
		where the square on the right is cartesian, $\sigma$ is given by quotienting by $S'_{\pi}$, and $\overline{p}$ is the induced quotient morphism. For all $h=1,\dots,d$, recall that $m_h$ is the number of blocks of $\pi$ of size exactly $h$. Then the morphism $Y_{\pi}\to Y_{\pi}/S_{\pi}$ in the bottom row of \eqref{descent} is an $\overline{S}_{\pi}$-Galois cover. The $\overline{S}_{\pi}$-action on $Y_{\pi}\times \prod_{j=1}^i \A^{|\pi_j|}$ is given by permutation of the factors corresponding to blocks of the same size. It follows by \'etale descent that the map $\overline{p} \colon \overline{p}^{-1}(Y_{\pi}/S_{\pi})\to Y_{\pi}/S_{\pi}$ is a vector bundle.  (A priori this vector bundle is only \'etale-locally trivial, but this implies that it is Zariski-locally trivial by Grothendieck's version of Hilbert Theorem 90.)
		
		Define $Z_{\pi} \coloneqq  \iota_{r-1}^{-1}(Y_{\pi})$ for every $\pi\in P(d)$, so that $\set{Z_{\pi}}_{\pi\in P(d)}$ is a locally closed stratification of $(\A^{r-1})^d_0\subset (\A^{r-1})^d$. Write $\sigma_{1,j}, \dots, \sigma_{|\pi_j|,j}$, $1\leq j\leq i$, for the standard basis of $\A^{|\pi_j|}$. By definition, $\sigma_{h,j}$ is the degree $h$ elementary symmetric polynomial on the coordinates of $(\A^1)^{|\pi_j|}$, which are the $a_{r,l}$ such that $l\in \pi_j$. The trivial vector bundle $p_0$ of \eqref{reduction} is the subbundle of $\iota_{r-1}^*p$ defined by the equation $a_{d,1}+\dots+a_{d,n}=0$. Over $Z_{\pi}$, this is exactly the inverse image of the subbundle of $Z_{\pi}\times (\prod_{j=1}^i \A^{|\pi_j|})$ given by the equation $\sigma_{1,1}+\sigma_{1,2}+\dots+\sigma_{1,j}=0$, which we call $Z_{\pi}\times (\prod_{j=1}^i \A^{|\pi_j|})_0$. Therefore, we have a cartesian square
		\[
		\begin{tikzcd}
			Z_{\pi}\times (\A^1)^d_0 \arrow[r, "\sigma"] \arrow[d, hook] & Z_{\pi}\times (\prod_{j=1}^i \A^{|\pi_j|})_0 \arrow[d, hook]  \\
			Z_{\pi}\times (\A^1)^d \arrow[r, "\sigma"] & Z_{\pi}\times (\prod_{j=1}^i \A^{|\pi_j|}). 
		\end{tikzcd}
		\]
  
		We now observe that $Z_{\pi}\times (\prod_{j=1}^i \A^{|\pi_j|})_0$ is a $\prod S_{m_h}$-equivariant vector subbundle of $Z_{\pi}\times (\prod_{j=1}^i \A^{|\pi_j|})$, and so descends to a subbundle of $\overline{p}^{-1}(Z_\pi/S_\pi)\to Z_\pi/S_\pi$. In other words,
		\[(Z_{\pi}\times (\A^1)^d_0)/S_{\pi}=\overline{p}^{-1}(Z_{\pi}/S_{\pi})\to Z_{\pi}/S_{\pi}\]
		is a vector bundle, as desired.
	\end{proof}
	
	Let $Y$ be a $k$-variety, $d$ be a positive integer, and let $S_d$ act on $(\A^1\times_k Y)^d$ by permutation of the $d$ factors $\A^1\times_k Y$. In \cite[Lemma 4.4]{G}, Totaro used a stratification argument to show that $[(\A^1 \times_k Y)^d/S_d]=[\A^d \times_k (Y^d/S_d)]$ in $K_0(\on{Var}_k)$. His arguments prove the following stronger statement.
	
		\begin{lemma}\cite[Lemma 4.4]{G}\label{totaro}
		Let $Y$ be a $k$-variety, $d$ be a positive integer, and let $S_d$ act on $(\A^1\times_k Y)^d$ and $Y^d$ by permutation of the $d$ factors $\A^1\times_k Y$ and $Y^d$, respectively. Then the projection $(\A^1\times_k Y)^d/S_d\to Y^d/S_d$ can be stratified by vector bundles (with height $1$).
	\end{lemma}

    We are ready to prove the main result of this subsection. In what follows, $S_d^m\rtimes S_m$ denotes the semidirect product, where $S_m$ acts on $S_d^m$ by permuting the factors.
	
	\begin{proposition}\label{zero-sum-bundle}
		Let $X$ be a quasi-projective $k$-variety, $d,m$ be non-negative integers, $S_m$ act on $X^m$ by permutation of factors, and $\on{pr}_i \colon X^m\to X$ be the $i$-th projection, for $i=1,\dots,m$. Let $E\to X$ be a vector bundle, and consider the $(S_d^m\rtimes S_m)$-equivariant vector bundle $\mathcal{E} \coloneqq \oplus_{i=1}^m \on{pr}_i^*((E^{\oplus d})_0)$ over $X^m$. Then the morphism $\mathcal{E}/(S_d^m\rtimes S_m)\to X^m/S_m$ can be stratified by vector bundles.	In particular, for every $x\in X^m/S_m$ we have $[(\mathcal{E}/(S_d^m\rtimes S_m))_{k(x)}]=0$ in $K_0(\on{Var}_{k(x)})/(\mathbb{L})$. 
	\end{proposition}
	
	\begin{proof}
		Define $\Sigma_{d,m} \coloneqq S_d^m\rtimes S_m$. For the purpose of proving that the morphism $\mathcal{E}/\Sigma_{d,m}\to X^m/S_m$ can be stratified by vector bundles, we are allowed to pass to a Zariski open cover of $X$. Indeed, let $x$ be a point of $X^{m}/S_{m}$, let $Z_x$ be the set-theoretic fiber of $X^{m}\to X^{m}/S_{m}$ at $x$, and set $Z'_x \coloneqq \cup_{i=1}^m\on{pr}_i(Z_x)\subset X$. Since $X$ is quasi-projective and $Z'_x$ is a finite subset of $X$, by \cite[00DS]{stacks} there exists an affine open subscheme $U\subset X$ such that $Z_x\subset U$, which implies that $Z_x$ is contained in $U^{m}$. Shrinking $U$ if necessary, we may also assume that $E|_U$ is trivial.  Therefore, we may suppose that $X=U$, that is, $E\cong \A^d_X$ is trivial, so that we have an isomorphism
		\[\mathcal{E}\cong ((\A^r)^d_0\times_kX)^{m}\]
		of $\Sigma_{d,m}$-equivariant vector bundles on $X^m$, where $r$ is the rank of $E$ and $(-)^m$ denotes $m$-fold fibered product over $\Spec(k)$.
		
		The projection maps $\A^r\to \A^{r-1}\to\dots\to \A^1\to \on{Spec}(k)$ given by forgetting the last component induce the maps
		\[((\A^r)^d_0\times_kX)^{m}/\Sigma_{d,m}\xrightarrow{p_r} ((\A^{r-1})^d_0\times_kX)^{m}/\Sigma_{d,m}\xrightarrow{p_{r-1}}\dots\to ((\A^{1})^d_0\times_kX)^{m}/\Sigma_{d,m}\xrightarrow{p_1} X^m/S_m.\]	
		In order to conclude, it suffices to show that $p_i$ can be stratified by vector bundles for all $i=1,\dots,d$. 
		
		By \Cref{zero-sum-over-k}, there exist a finite set $J$ and an $S_d$-equivariant locally closed stratification $\set{U_j}_{j\in J}$ of $(\A^{i-1})^d_0/S_d$ such that the restriction of $((\A^{i})^d_0)/S_d\to ((\A^{i-1})^d_0)/S_d$ to each $U_j$ is a trivial vector bundle. Let $X_j \coloneqq  U_j\times X$. Since $S_d$ acts trivially on $X$, the restriction of $((\A^{i})^d_0\times_kX)/S_d\to ((\A^{i-1})^d_0\times_kX)/S_d$ to each $X_j$ is a trivial vector bundle. It follows that $\set{X_{j_1}\times \dots\times X_{j_m}}_{(j_1,\dots,j_m)\in J^m}$ is a $\Sigma_{d,m}$-equivariant stratification of $((\A^{i-1})^d_0\times X)^m$ such that the restriction of $((\A^{i})^d_0\times_kX)^{m}/S_d^m\to ((\A^{i-1})^d_0\times_kX)^{m}/S_d^m$ to each stratum is a trivial vector bundle:
		\[\A^{d-1}_{X_{j_1}}\times\dots\times \A^{d-1}_{X_{j_m}}\to X_{j_1}\times \dots\times X_{j_m}.\] 
  
		It remains to take the quotient by $S_m$. For every $\mathbf{j}=(j_1,\dots,j_m)\in J^m$, let $\pi(\mathbf{j})$ be the partition of $\set{1,\dots,m}$ such that $s$ and $t$ are in the same block of $\pi$ if and only if $j_s=j_t$. The symmetric group $S_m$ acts on the stratification $\set{X_{j_1}\times \dots\times X_{j_m}}_{\mathbf{j}\in J^m}$ via its component-wise permutation action on $J^m$. The stabilizer of $X_{j_1}\times \dots\times X_{j_m}$ is the direct product of the symmetric groups of each block of $\pi(\mathbf{j})$. By \Cref{reduce-to-one-component}, it is enough to show that each
		\[(\A^{d-1}_{X_{j_1}}\times\dots\times \A^{d-1}_{X_{j_m}})/S'_{\pi(\mathbf{j})}\to (X_{j_1}\times \dots\times X_{j_m})/S'_{\pi(\mathbf{j})}\] 
		can be stratified by vector bundles. The latter map factors as the composition
		\[(\A^{d-1}_{X_{j_1}}\times\dots\times \A^{d-1}_{X_{j_m}})/S'_{\pi(\mathbf{j})}\to (\A^{d-2}_{X_{j_1}}\times\dots\times \A^{d-2}_{X_{j_m}})/S'_{\pi(\mathbf{j})}\to\dots\to (X_{j_1}\times \dots\times X_{j_m})/S'_{\pi(\mathbf{j})},\]
		where each map forgets the last component of each of the $\A^h_{X_{j_s}}$.
		Following the decomposition of $\set{1,\dots,m}$ into blocks of $\pi(\mathbf{j})$, each of these maps is a product of maps of the form
		\[((\A^h_{X_{j_s}})^{m(j_s)})/S_{m(j_s)}\to ((\A^{h-1}_{X_{j_s}})^{m(j_s)})/S_{m(j_s)}.\]
		Here $m(j_s)$ is the number of $t\in\set{1,\dots,m}$ such that $j_s=j_t$, the notation $(-)^{m(j_s)}$ indicates  $m(j_s)$-fold direct product over $\Spec(k)$, and $S_{m(j_s)}$ acts by permuting the factors $\A^h_{X_{j_s}}$. Setting $Y=\A^{h-1}_{X_{j_s}}$ and $m=m(j_s)$ puts this map in the form of \Cref{totaro}, which gives the conclusion.
	\end{proof}

	\subsection{The Polydiagonal Compactification}
    \label{subsection-polydiag}
	
	Let $X$ be a smooth quasi-projective variety over a field $k$ of characteristic zero, and let $n \geq 1$ be an integer. The symmetric power $\on{Sym}^n(X)$ is given by the quotient $X^n/S_n$, which does not act with diagonalizable stabilizers for $n \geq 3$.  However, there is a canonical equivariant blowup of $X^n$ such that the action does have diagonalizable stabilizers; this is the {\it polydiagonal compactification} of configuration space due to Ulyanov.  We'll review some facts about this construction below and refer to \cite{Ulyanov} for more details.
	
    We'll continue to use the notation for partitions from the beginning of this section.  For any partition $\pi \in P(n)$, define the {\it polydiagonal} $\Delta^{\pi} \subset X^n$ to be the closed subscheme parametrizing $n$-tuples $(x_1,\dots,x_n)$ such that $x_k = x_l$ whenever $k$ and $l$ belong to the same block of the partition $\pi$.  Note that $\Delta^{\pi} \cong X^i$, where $i$ is the number of blocks of $\pi$.  
	
	The polydiagonal compactification $X\langle n \rangle$ is constructed by blowing up $X^n$ as follows:
	
	\begin{definition}{\cite[Definition-Theorem 2.1]{Ulyanov}}
		The polydiagonal compactification $X\langle n \rangle$ of the configuration space $F(X,n)$ of ordered points is a sequence of $n-1$ blowups of $X^n$:
		$$X \langle n \rangle = Y_{n-1} \xrightarrow{\alpha_{n-1}} Y_{n-2} \xrightarrow{\alpha_{n-2}} \cdots \xrightarrow{\alpha_2} Y_1 \xrightarrow{\alpha_1} Y_0 = X^n.$$
		Here $Y_i \rightarrow Y_{i-1}$ is the blowup of the disjoint union of proper transforms $\Delta^{\pi}_{i-1}$ of polydiagonals $\Delta^{\pi}$ in $X^n$ corresponding to partitions $\pi$ with exactly $i$ blocks.
	\end{definition}
	
		Note that the polydiagonals $\Delta^{\pi}$ for $\pi$ of a particular type with exactly $i\geq 1$ blocks are not disjoint in $X^n$, but their proper transforms become disjoint by the $i-1$ step.
		
		The action of $S_n$ on $X \langle n \rangle$ has abelian isotropy subgroups by \cite[Theorem 3.11]{Ulyanov}.  In fact, the proof shows even more: the isotropy subgroup is contained in a direct sum of copies of $k^{\times}$, hence is diagonalizable. 

  We are now ready to prove \Cref{compactification-theorem}.
	
	\begin{proof}[Proof of \Cref{compactification-theorem}]
		We'll break the morphism $X \langle n \rangle /S_n \rightarrow X^n/S_n$ into a composition of $Y_i/S_n \rightarrow Y_{i-1}/S_n$ for all $i=1,\dots,n-1$.  We can subdivide this process further and consider the blowup of the disjoint union of all partitions of a particular type one at a time.  Choose a type $\mathbf{a}$ with $i$ blocks, $1 \leq i < n$; then $\mathbb{P}(\mathcal{N}_{Z_{\mathbf{a}}/Y_{i-1}})$ is the exceptional divisor in $Y_i$ obtained by blowing up the disjoint union $Z_{\mathbf{a}}$ of all proper transforms $\Delta^{\pi}_{i-1}$ in $Y_{i-1}$ of polydiagonals of type $\mathbf{a}$. Then both $Z_{\mathbf{a}}$ and $\mathbb{P}(\mathcal{N}_{Z_{\mathbf{a}}/Y_{i-1}})$ are $S_n$-invariant.
		
		By \Cref{closure}, the proof will be complete if we can show that the morphism \[\mathbb{P}(\mathcal{N}_{Z_{\mathbf{a}}/Y_{i-1}})/S_n \rightarrow Z_{\mathbf{a}}/S_n\] has $\mathbb{L}$-rational fibers.  Using \Cref{reduce-to-one-component}, we can further reduce to showing that \[\mathbb{P}(\mathcal{N}_{\Delta^{\pi}_{i-1}/Y_{i-1}})/S_\pi \to \Delta^{\pi}_{i-1}/S_\pi\] has $\mathbb{L}$-rational fibers for each partition $\pi$. By the blowup formula for the normal bundle \cite[Proposition B.6.10]{Fulton}, there exists an $S_{\pi}$-equivariant line bundle $L_{\pi,i-1}$ on $\Delta^{\pi}_{i-1}$ and an $S_{\pi}$-equivariant isomorphism 
		\[\mathcal{N}_{\Delta^{\pi}_{i-1}/Y_{i-1}}\cong (\alpha_1\circ\dots\circ \alpha_{i-1})^*(\mathcal{N}_{\Delta^{\pi}/X^n})\otimes L_{\pi,i-1}.\] 
		Therefore, the composition $\alpha_1\circ\dots\circ \alpha_{i-1}$ induces a commutative $S_{\pi}$-equivariant square 
		\[
		\begin{tikzcd}
			\mathbb{P}(\mathcal{N}_{\Delta^{\pi}_{i-1}/Y_{i-1}}) \arrow[r] \arrow[d] & \mathbb{P}(\mathcal{N}_{\Delta^{\pi}/X^n}) \arrow[d] \\ 
			\Delta^{\pi}_{i-1} \arrow[r] & 	\Delta^{\pi}.
		\end{tikzcd}
		\]
		The fiber of a pullback bundle over a point is isomorphic to the base change of the fiber over the image point, and the projectivization (with its $S_{\pi}$-action) is insensitive to twisting by $L_{\pi,i-1}$.  Therefore, it suffices to show instead that $\mathbb{P}(\mathcal{N}_{\Delta^{\pi}/X^n})/S_\pi \to \Delta^{\pi}/S_\pi$ has $\mathbb{L}$-rational fibers.
		
		Since $k$ is of characteristic zero, by \cite[Proposition 1.9]{GIT}, finite group quotients commute with arbitrary base change. Moreover, projectivization commutes with arbitrary base change. Choose a point $x$ of $\Delta^{\pi}/S_\pi$ and let $W_x\subset \Delta^{\pi}$ be the fiber of $x$, so that $W_x/S_\pi=\on{Spec}(k(x))$. We deduce that
		\begin{equation}\label{quot-commutes}(\mathbb{P}(\mathcal{N}_{\Delta^{\pi}/X^n})/S_\pi)_{k(x)}\cong (\mathbb{P}(\mathcal{N}_{\Delta^{\pi}/X^n})_{W_x})/S_\pi\cong \mathbb{P}((\mathcal{N}_{\Delta^{\pi}/X^n})_{W_x})/S_\pi.\end{equation}
  
		For the first isomorphism, we have applied \cite[Proposition 1.9]{GIT} to $\mathbb{P}(\mathcal{N}_{\Delta^{\pi}/X^n})\to \Delta^{\pi}/S_{\pi}$. 
    For the second isomorphism, we have used the commutativity of projectivization and base change.
		
		From (\ref{quot-commutes}) and \Cref{cor-tautological-bundle}, we deduce that $[(\mathbb{P}(\mathcal{N}_{\Delta^{\pi}/X^n})/S_\pi)_{k(x)}]=1$ in $K_0(\on{Var}_{k(x)})/(\mathbb{L})$ is equivalent to $[((\mathcal{N}_{\Delta^{\pi}/X^n})_{W_x})/S_\pi]=0$ in $K_0(\on{Var}_{k(x)})/(\mathbb{L})$. In order to conclude, it now suffices to show that $[(\mathcal{N}_{\Delta^{\pi}/X^n}/S_\pi)_{k(x)}] = 0$ in $K_0(\on{Var}_{k(x)})/(\mathbb{L})$.
		
		The closed embedding $\Delta^{\pi}\subset X^n$ determines an $S_\pi$-equivariant short exact sequence of vector bundles on $\Delta^{\pi}$:
		\[0\to \oplus_{j=1}^ip_j^*\mathcal{T}_X\xrightarrow{\iota_{\pi}} \oplus_{l=1}^n q_l^*\mathcal{T}_X|_{\Delta^{\pi}}\to \mathcal{N}_{\Delta^{\pi}/X^n}\to 0.\]
		(Recall that $i$ is the number of blocks of $\pi$.) Here $p_j \colon \Delta^{\pi}\to X$ and $q_l \colon X^n\to X$ denote the $j$-th and $l$-th projection, respectively; note that $q_l|_{\Delta^{\pi}}=p_j$ if and only if $l\in \pi_j$. The map $\iota_{\pi}$ sends $(v_1,\dots,v_i)$ to $(w_1,\dots,w_n)$, where $w_l=v_j$ whenever $l\in \pi_j$. This map has a retraction given by the $S_{\pi}$-equivariant map $\rho_{\pi}$ which sends $(w_1,\dots,w_n)$ to 
		$(v_1,\dots,v_i)$, where $v_j=\frac{1}{|\pi_j|}\sum_{l\in \pi_j}w_l$ for all $j=1,\dots,i$. We obtain an $S_{\pi}$-equivariant isomorphism of vector bundles
		\[\mathcal{N}_{\Delta^{\pi}/X^n}\cong \on{coker}\iota_{\pi}\cong \ker \rho_{\pi}\cong \oplus_{j=1}^i p_j^*(\mathcal{T}_X^{\oplus (|\pi_j|-1)}),\]
		where on the right side, each $S_{\pi_j}$ acts on the corresponding summand $\mathcal{T}_X^{\oplus (|\pi_j|-1)}$ as $(\mathcal{T}_X^{\oplus (|\pi_j|)})_0$, that is, as the kernel of the addition map $\mathcal{T}_X^{\oplus |\pi_j|}\to \mathcal{T}_X$, where $\mathcal{T}_X^{\oplus |\pi_j|}$ is the standard permutation representation of $S_{\pi_s}$, and trivially on the other summands, and each $S_{m_h}$ acts by permuting the factors corresponding to the blocks $\pi_j$ which satisfy $|\pi_j|=h$. (Recall that $m_h$ is defined as the number of blocks of $\pi$ of size exactly $h$.) We can write the map in question \[\mathcal{N}_{\Delta^{\pi}/X^n}/S_\pi=\oplus_{j=1}^i p_j^*(\mathcal{T}_X^{\oplus (|\pi_j|-1)})/S_\pi\to \Delta^{\pi}/S_\pi\] as a direct product of maps of the form 
		\[\left(\oplus_{j \colon |\pi_j|=h} \ p_j^*(\mathcal{T}_X^{h})_0 \right)/(S_{h}^{m_h}\rtimes S_{m_h})\to X^{m_h}/S_{m_h},\]
		where $h\in \set{1,\dots,n}$ satisfies $m_h\geq 1$. (The condition $m_h\geq 1$ appears because we considered the groups $S_m$ for $m\geq 1$ only.) The conclusion then follows from \Cref{zero-sum-bundle}.	
	\end{proof}
	
	By \Cref{singcompclosure}, to answer \Cref{mainq}(i) it would therefore be enough to determine if $X \langle n \rangle / S_n$ has $\mathbb{L}$-rational singularities.

\end{document}